\documentclass[11pt]{article} 

\usepackage[utf8]{inputenc} 

\usepackage{geometry} 
\usepackage{graphicx} 

\usepackage{booktabs} 
\usepackage{array} 
\usepackage{paralist} 
\usepackage{verbatim} 
\usepackage{subfig} 
\usepackage{amssymb} 
\usepackage{nccmath} 
\usepackage{ntheorem} 
\usepackage{latexsym}
\usepackage{tikz}
\usepackage{mathrsfs}

\usepackage{fancyhdr} 
\usepackage{sectsty}
\allsectionsfont{\sffamily\mdseries\upshape} 

\usepackage[nottoc,notlof,notlot]{tocbibind} 
\usepackage[titles,subfigure]{tocloft} 

\newtheorem{thm}{Theorem}[section]

\newtheorem{cor}[thm]{Corollary}
\newtheorem{lem}[thm]{Lemma}

\newtheorem{conj}[thm]{Conjecture}
\theorembodyfont{\normalfont}

\newtheorem*{ack}{Acknowledgements.}
\theoremheaderfont{\it}
\newtheorem*{prf}{Proof.}
\newtheorem*{mainprf}{Proof of Theorem \ref{mainres}.}
\theoremheaderfont{\slshape}

\newtheorem{rem}[thm]{Remark}

\def\qed{$\hfill\Box$}

\begin{document}
\begin{center}
{\large \textbf{The action of full twist on the superpolynomial for torus knots}}\\
\vspace{20pt}
{\large Keita Nakagane}\footnote[2]{Department of Mathematics, Tokyo Institute of Technology. E-mail: nakagane.k.aa@m.titech.ac.jp}
\end{center}
\vspace{10pt}
\noindent
\textbf{Abstract.}
We show, using Mellit's recent results, that K\'alm\'an's full twist formula for the HOMFLY polynomial can be generalized to a formula for superpolynomials in the case of positive toric braids.

\section{Introduction}
The HOMFLY polynomial \cite{Jo} $P(K)(Q,\alpha)$ is an oriented link invariant, which is determined by the skein relation
\[
\alpha^{-1}P(
\begin{tikzpicture}[baseline=1]
\draw[->] (0,0)--(0.3,0.3);
\draw (0.3,0)--(0.2,0.1);
\draw[->] (0.1,0.2)--(0,0.3);
\end{tikzpicture}
)-\alpha P(
\begin{tikzpicture}[baseline=1]
\draw[->] (0.3,0)--(0,0.3);
\draw (0,0)--(0.1,0.1);
\draw[->] (0.2,0.2)--(0.3,0.3);
\end{tikzpicture}
)=(Q^{-1}-Q)P(
\begin{tikzpicture}[baseline=1]
\draw [bend right = 45,->] (0,0) to (0,0.3);
\draw [bend left = 45,->] (0.3,0) to (0.3,0.3);
\end{tikzpicture}
)
\]
and $P(\mathrm{trivial\ knot})=1$.
This invariant includes the Alexander polynomial $P(K)(Q,1)$, and it also includes the $\mathfrak{sl}(N)$ polynomials $P(K)(Q,Q^N)$ which are commonly referred to as quantum invariants.
Here the $\mathfrak{sl}(2)$ polynomial is the Jones polynomial.

The Khovanov--Rozansky homologies \cite{KR1, KR2} are homological link invariants which categorify the HOMFLY polynomial and the $\mathfrak{sl}(N)$ polynomials.
This means that the homologies are multi-graded vector spaces depending only on the link type, and that the decategorified polynomial invariants can be obtained as certain graded Euler characteristics of the homologies.  
The HOMFLY homology $H(K)$ is triply graded, whereas the $\mathfrak{sl}(N)$ homologies $H_N(K)$ are doubly graded.

The superpolynomial $\mathcal{P}(K)(Q,\alpha,T)$ is a conjectural knot invariant, introduced by Dunfield, Gukov and Rasmussen \cite{DGR}, satisfying $\mathcal{P}(K)(Q,\alpha,-1) = P(K)(Q,\alpha)$.
The polynomial $\mathcal{P}(K)$ is the Poincar\'e polynomial of a conjectural triply graded vector space $\mathcal{H}(K)$, which has differentials $\{d_N\colon \mathcal{H}(K) \to \mathcal{H}(K)\}_{N \in \mathbb{Z}}$ satisfying some properties.
For example, the homology $H(\mathcal{H}(K),d_N)$ should be isomorphic to $H_N(K)$ for all $N>0$.
Furthermore, $H(\mathcal{H}(K),d_0)$ should be isomorphic to the knot Floer homology of $K$.
It has been conjectured that  the HOMFLY homology $H(K)$ can be taken as $\mathcal{H}(K)$. 
The differential $d_N$ for $N>0$ on $H(K)$ was constructed by Rasmussen \cite{Ra}, as the first differential of a spectral sequence from $H(K)$ to $H_N(K)$, and the property stated above was proved by Naisse and Vaz \cite{NV}. Although not all properties are proved yet, we call the Poincar\'e polynomial of $H(K)$ the superpolynomial of $K$ here.

Our main interests are in the action of full twists on braids and link invariants of closures.
K\'alm\'an \cite{Ka1} showed the following full twist formula for the HOMFLY polynomial:
\begin{description}
\item[] the coefficient of $\alpha^{e-n+1}$ in $P(\widehat{\beta})$ $=$ $(-1)^{n-1}$ $\times$ the coefficient of $\alpha^{e+n-1}$ in $P(\widehat{\beta\Delta_{n}^2})$,
\end{description}
where $\beta$ is a braid on $n$ strands, $e$ is the exponent sum of $\beta$ and $\Delta_n$ is the positive half twist of $n$ strands. 
We conjecture that this formula can be generalized to the superpolynomial:
\begin{description}
\item[] the coefficient of $\alpha^{e-n+1}$ in $\mathcal{P}(\widehat{\beta})$ $=$ $T^{(n^2-1)}$ $\times$ the coefficient of $\alpha^{e+n-1}$ in $\mathcal{P}(\widehat{\beta\Delta^2})$.
\end{description}
We will show in this paper (Theorem \ref{mainres}) that this is true if $\beta$ is a positive toric braid by using Mellit's result \cite{Me}.
His result was first conjectured by Gorsky and Negut \cite{GN} via physical considerations. 

The paper is organized as follows.
In section 2, we explain our notation.
In section 3, we recall K\'alm\'an's full twist formula and state our conjecture more precisely.
We also prove our conjecture in the case of positive toric braids.
\begin{ack}
I am sincerely grateful to my supervisor, associate professor Tam\'as K\'alm\'an for his encouragement and valuable advice. 
\end{ack}
\section{Notation}
We will use nearly the same notation as Elias and Hogancamp \cite{EH}.
By \textit{superpolynomial}, we mean the reduced superpolynomial in Appendix A of \cite{EH}, which agrees with what we described in the introduction.
The superpolynomial of a link $L$, denoted by $\mathcal{P}(L)(Q,\alpha,T)$, is a Laurent series in $Q$, $\alpha$ and $T$ satisfying $\mathcal{P}(L)(Q,\alpha,-1) = P(L)(Q,\alpha)$.

Next we recall the Morton--Franks--Williams inequality \cite{FW, Mo}.
For a link $L$, let $d_{-}(L)$ (resp.\ $d_{+}(L)$) be the lowest (resp.\ highest) degree of $\alpha$ in $P(L)$.
\begin{thm}[\cite{FW, Mo}]
Let $L$ be a link. If $\beta$ is a braid whose closure is $L$ then
\[
e(\beta)-n(\beta)+1 \ \le \ d_{-}(L) \ \le \ d_{+} (L) \ \le \ e(\beta)+n(\beta)-1,
\]
where $e(\beta)$ is the exponent sum of $\beta$ and $n(\beta)$ is the number of strands of $\beta$.
\end{thm}
Let us also mention that essentially the same degree bounds also apply to the superpolynomial.
For a braid $\beta$, the coefficient of $\alpha^{e(\beta)-n(\beta)+1}$ (resp.\ $\alpha^{e(\beta)+n(\beta)-1}$) in $P(\widehat{\beta})$ will be denoted by $P_{-}(\beta)$ (resp.\ $P_{+}(\beta)$).
For the superpolynomial, $\mathcal{P}_{-}(\beta)$ and $\mathcal{P}_{+}(\beta)$ are similarly defined.

\section{Adding a full twist}
K\'alm\'an \cite{Ka1} showed the \textit{full twist formula} for the HOMFLY polynomial. 
\begin{thm}[\cite{Ka1}]
For an $n$-strand braid $\beta$, we have $P_{-}(\beta) = (-1)^{n-1}P_{+}(\beta\Delta_n^{2})$, where $\Delta_n$ is the positive half twist of $n$ strands.
\end{thm}
It is reasonable to consider the generalization of this formula to the superpolynomial.
Now we have the following conjecture.
\begin{conj}\label{conj}
For an $n$-strand braid $\beta$, we have $\mathcal{P}_{-}(\beta) = T^{(n^2-1)}\mathcal{P}_{+}(\beta\Delta_n^{2})$.
\end{conj}
We do not have a proof of the conjecture yet, but we have one result supporting it.
\begin{thm}\label{mainres}
Let $(m,n)$ be a pair of coprime positive integers.
For the $(m,n)$-toric braid $\tau_{m,n} = (\sigma_1\cdots \sigma_{n-1})^m$ on $n$ strands, we have 
\[
\mathcal{P}_{-}(\tau_{m.n})\,=\,T^{(n^2-1)}\mathcal{P}_{+}(\tau_{m,n}\Delta^{2}_n)\,=\,T^{(n^2-1)}\mathcal{P}_{+}(\tau_{m+n,n}).
\]
\end{thm}

The proof relies on Mellit's result for positive torus knots \cite{Me}, which was conjectured by Gorsky and Negut \cite{GN} and by Oblomkov, Rasmussen and Shende \cite{ORS}.
Let us recall the formula.

Let $(m,n)$ be a pair of coprime positive integers. 
An \textit{$(m,n)$-Dyck path} is a lattice path $\gamma$ in $\mathbb{R}^2$ from $(0,0)$ to $(m,n)$ obtained by concatenating $m$ \textit{horizontal steps} $(1,0)$ and $n$ \textit{vertical steps} $(0,1)$ in any order so that $\gamma$ is sitting above the \textit{diagonal line} $y =(n/m)x$.
Let $D_{m,n}$ be the set of $(m,n)$-Dyck paths.
By \textit{parallel line} we mean a line in $\mathbb{R}^2$ parallel to the diagonal line.
For a point $p \in \mathbb{R}^2$, we denote by $l(p)$ the parallel line passing through $p$.

Let $\gamma \in D_{m,n}$ be a Dyck path.
We denote by $\mathrm{area}(\gamma)$ the number of complete $1 \times 1$ lattice squares between $\gamma$ and the diagonal line.
We define the sets $O(\gamma)$ and $H(\gamma)$ by
\begin{align*}
&O(\gamma) = \left\{ (r_h,r_v) \middle|
\begin{array}{l}
r_h \text{ (resp.\ } r_v\text{) is a horizontal (resp.\ vertical) step of }\gamma\\
r_v \text{ appears after }r_h\text{ in }\gamma
\end{array}
\right\}\\
&H(\gamma) = \left\{ (r_h,r_v) \in O(\gamma) \middle|
\begin{array}{l}
\text{there exists a parallel line }l\text{ intersecting both }r_h\text{ and }r_v
\end{array}
\right\}
\end{align*}
and let $h(\gamma) = |H(\gamma)|$ be the number of pairs in $H(\gamma)$.
An \textit{outer vertex} of $\gamma$ is a lattice point in $\gamma$ which is just after a vertical step and just before a horizontal step in $\gamma$.
For an outer vertex $p$, we denote by $k(p)$ the number of horizontal steps (or vertical steps) in $\gamma$ intersecting $l(p)$.
(The line $l(p)$ touches two steps at $p$ but we do not count them.)
We remark that $k(p_0) = 0$ where $p_0$ is the most distant outer vertex of $\gamma$ from the diagonal line. 
Let $V(\gamma)$ be the set of all outer vertices of $\gamma$ except for $p_0$.
It is obvious that $|V(\gamma)| \le n-1$.
We will say that $\gamma$ is \textit{rugged} if $|V(\gamma)| = n-1$.
Let $D^{*}_{m,n}$ be the set of rugged $(m,n)$-Dyck paths.

An example of a Dyck path for $(m,n) = (5,4)$ is shown below.
For this Dyck path $\gamma$, we have $\mathrm{area}(\gamma)=2$, $h(\gamma) = 4$ and $V(\gamma) = \{(0,2), (3,4)\}$.
(The outer vertex (1,3) is the most distant vertex from the diagonal line $y = (4/5)x$.) 
Moreover, $k(0,2) = 1$ and $k(3,4) = 2$.
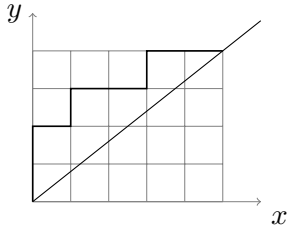
\begin{figure}[h]
\centering
\begin{tikzpicture}
\draw[help lines, <->] (0,2.5)--(0,0)--(3,0);
\node[below right] at (3,0) {$x$};
\node[left] at (0,2.5) {$y$};
\draw[help lines, step=0.5] (0,0) grid (2.5,2);
\draw (0,0)--(3,2.4);
\draw[semithick] (0,0)--(0,1)--(0.5,1)--(0.5,1.5)--(1.5,1.5)--(1.5,2)--(2.5,2);
\end{tikzpicture}
\caption{A $(5,4)$-Dyck path $\gamma$}
\label{fig1}
\end{figure}

Now we are in a position to state Mellit's result.
\begin{thm}[\cite{Me}]
Let $(m,n)$ be a pair of coprime positive integers.
The superpolynomial of the $(m,n)$-torus knot $T_{m,n}$ is given by 
\[
\mathcal{P}(T_{m,n}) = (T^{-1}\alpha)^{(m-1)(n-1)}\sum_{\gamma \in D_{m,n}}q^{\mathrm{area}(\gamma)}t^{h(\gamma)}\prod_{p \in V(\gamma)}(1+T^{-1}\alpha^2t^{-k(p)}),
\]
where $q = Q^2$ and $t = T^2Q^{-2}$.
\end{thm}
\begin{rem}
Let us explain the difference of notation with \cite{Me}.
If we denote by $q'$, $a'$ and $t'$ the variables $q$, $a$ and $t$ in \cite{Me} respectively then 
\[
q' = t, \qquad a'=-T^{-1}\alpha^2, \qquad t'=q.
\]
$\mathcal{P}$ in \cite{Me} is the product of our $\mathcal{P}$ and $(-1)^{e(\beta)-n(\beta)+\mu(\beta)}(\alpha Q^{-1})^{\mu(\beta)}/(1-Q^2)$, where $\mu(\beta)$ is the number of components of the closure $\widehat{\beta}$.
\end{rem}
\begin{cor}\label{cor}
Let $(m,n)$ be a pair of coprime positive integers. Then
\begin{align*}
\mathcal{P}_{-}(\tau_{m,n}) &= T^{-(m-1)(n-1)}\sum_{\gamma \in D_{m,n}}q^{\mathrm{area}(\gamma)}t^{h(\gamma)}\\
\mathcal{P}_{+}(\tau_{m,n}) &= T^{-m(n-1)}\sum_{\gamma \in D^{*}_{m,n}}q^{\mathrm{area}(\gamma)}t^{h(\gamma)-\sum_{p \in V(\gamma)}k(p)}\vphantom{a}_.
\end{align*}
\end{cor}
\begin{prf}
The first equation is obvious from $e(\tau_{m,n})-n(\tau_{m,n})+1 = (m-1)(n-1)$. 
On the other hand, $e(\tau_{m,n})+n(\tau_{m,n})-1 = (m-1)(n-1)+2(n-1)$ hence the contributions to $\mathcal{P}_{+}(\tau_{m,n})$ of non-rugged paths are $0$ and those of rugged paths are as shown in the second equation. \qed
\end{prf}
For the proof of our result, we prepare the following lemmas.
\begin{lem}\label{LEM1}
Let $(m,n)$ be a pair of coprime positive integers and let $\gamma \in D_{m,n}$ be a Dyck path.
Then
\[
\mathrm{area}(\gamma) = \frac{(m-1)(n-1)}{2}-|O(\gamma)|.
\]
\end{lem}
\begin{prf}
Let $\gamma_0 \in D_{m,n}$ be the Dyck path defined by concatenating $n$ vertical steps first and $m$ horizontal steps after them.
It is easily verified that $\mathrm{area}(\gamma_0) = (m-1)(n-1)/2$.
By the definition of $\gamma_0$, the number of lattice squares between $\gamma_0$ and $\gamma$ is $\{\mathrm{area}(\gamma_0)-\mathrm{area}(\gamma)\}$.
On the other hand, those squares are in a one-to-one correspondence with pairs in $O(\gamma)$; for each such square, there is a horizontal step $r_h$ (resp.\ vertical step $r_v$) in $\gamma$ parallel to the bottom edge (resp.\ left edge) of the square such that $r_v$ appears after $r_h$ in $\gamma$.
Therefore, $\mathrm{area}(\gamma) = \mathrm{area}(\gamma_0)-|O(\gamma)| = (m-1)(n-1)/2-|O(\gamma)|$. \qed
\end{prf}
\begin{lem}\label{LEM2}
Let $(m,n)$ be a pair of coprime positive integers and let $\gamma \in D_{m,n}$ be a Dyck path. 
Let $(r_h,r_v) \in O(\gamma)$ and suppose that $r_h$ is a step to $(a,b)$ and $r_v$ is a step from $(a',b')$.
Then the following claims {\rm(1A)} and {\rm(1B)}, as well as {\rm(2A)} and {\rm(2B)}, are equivalent:
\begin{fleqn}
\begin{description}
\item[(1A)] $l(a-1,b)$ intersects $r_v$
\item[(1B)] $-n < n(a'-a)-m(b'-b) < m-n$
\item[(2A)] $l(a',b'+1)$ intersects $r_h$
\item[(2B)] $m-n < n(a'-a)-m(b'-b) < m$.
\end{description}
In particular, as {\rm(1B)} and {\rm(2B)} are contradictory, {\rm(1A)} and {\rm(2A)} do not hold at the same time.
Moreover, $(r_h,r_v) \in H(\gamma)$ if and only if {\rm(1A)} or {\rm(2A)} holds.
\end{fleqn} 
\end{lem}
\begin{prf}
It is easy to verify the equivalence statements. 
The ``only if'' part of the last statement can be proved by lifting up a line intersecting $r_h$ and $r_v$ until it passes through $(a-1,b)$ or $(a',b'+1)$. \qed
\end{prf}
By lemma \ref{LEM2}, $H(\gamma)$ can be partitioned into two sets $H_1(\gamma)$ and $H_2(\gamma)$.
To be precise, $H_1(\gamma)$ (resp.\ $H_2(\gamma)$) is the set of pairs in $H(\gamma)$ satisfying condition (1A) (resp.\ (2A)).
 
Now we prove our main result.
\begin{mainprf}
By Corollary \ref{cor},
\[
\mathcal{P}_{-}(\tau_{m,n}) = T^{-(m-1)(n-1)}\sum_{\gamma \in D_{m,n}}q^{\mathrm{area}(\gamma)}t^{h(\gamma)}
\]
and
\[
T^{(n^2-1)}\mathcal{P}_{+}(\tau_{m+n,n}) = T^{-(m-1)(n-1)}\sum_{\gamma \in D^{*}_{m+n,n}}q^{\mathrm{area}(\gamma)}t^{h(\gamma)-\sum_{p \in V(\gamma^{*})}k(p)}\vphantom{a}_.
\]
We have a natural bijection between $D_{m,n}$ and $D^{*}_{m+n,n}$: a path $\gamma\in D_{m,n}$ corresponds to the path $\gamma^{*}\in D^{*}_{m+n,n}$ obtained by inserting one horizontal step just after each of the vertical steps of $\gamma$.
(We show an example below, namely $\gamma^{*}$ for the Dyck path $\gamma$ in Figure \ref{fig1}.
The thickened horizontal steps are the newly inserted steps of $\gamma^{*}$.)  
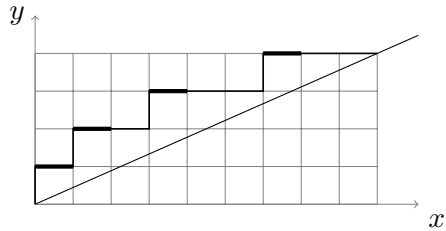
\begin{figure}[h]
\centering
\begin{tikzpicture}
\draw[help lines, <->] (0,2.5)--(0,0)--(5.04,0);
\node[below right] at (5.04,0) {$x$};
\node[left] at (0,2.5) {$y$};
\draw[help lines, step=0.5] (0,0) grid (4.5,2);
\draw (0,0)--(5.04,2.24);
\draw[semithick] (0,0)--(0,0.5)--(0.5,0.5)--(0.5,1)--(1.5,1)--(1.5,1.5)--(3,1.5)--(3,2)--(4.5,2);
\draw[ultra thick] (0,0.5)--(0.5,0.5);
\draw[ultra thick] (0.5,1)--(1,1);
\draw[ultra thick] (1.5,1.5)--(2,1.5);
\draw[ultra thick] (3,2)--(3.5,2);
\end{tikzpicture}
\caption{$\gamma^{*}$ for $\gamma$ in Figure \ref{fig1}}
\label{fig2}
\end{figure}

A horizontal step $r_h$ to $(a,b)$ in $\gamma$ becomes a horizontal step $r_h^{*}$ to $(a+b,b)$ in $\gamma^{*}$, and a vertical step $r_v$ from $(a',b')$ in $\gamma$ becomes a vertical step $r_v^{*}$ from $(a'+b',b')$ in $\gamma^{*}$.
This correspondence between $D_{m,n}$ and $D^{*}_{m+n,n}$ is truly a bijection since $nx-my<0$ is equivalent to $n(x+y)-(m+n)y<0$.
Consequently, it is sufficient to prove
\begin{description}
\item{(A)} $\mathrm{area}(\gamma) = \mathrm{area}(\gamma^{*})$
\item{(B)} $\displaystyle h(\gamma) = h(\gamma^{*})-\sum_{p \in V(\gamma^{*})}k(p)$
\end{description}
for each $\gamma \in D_{m,n}$.

Let us prove the equation (A) first. 
By Lemma \ref{LEM1}, 
\begin{align*}
\mathrm{area}(\gamma) &= \frac{(m-1)(n-1)}{2}-|O(\gamma)|\\
\mathrm{area}(\gamma^{*}) &= \frac{(m+n-1)(n-1)}{2}-|O(\gamma^{*})|.
\end{align*}
On the other hand, $|O(\gamma^{*})| = |O(\gamma)|+n(n-1)/2$ by the construction of $\gamma^{*}$ and thus $\mathrm{area}(\gamma) = \mathrm{area}(\gamma^{*})$.

Now we prove the equation (B).
Consider the injection $\iota: O(\gamma) \to O(\gamma^{*})$ induced by the correspondence $*$ and identify $(r_h,r_v) \in O(\gamma)$ with $\iota(r_h,r_v) = (r_h^{*},r_v^{*}) \in O(\gamma^{*})$.
Since 
\[
-n < n(a'-a)-m(b'-b) < m
\] 
is equivalent to
\[
-n < n\{(a'+b')-(a+b)\}-(m+n)(b'-b) < (m+n)-n,
\]
we have $H(\gamma) = H_{1}(\gamma^{*}) \cap O(\gamma)$ by Lemma \ref{LEM2}.
Hence it is sufficient to prove that
\begin{equation}\label{ess}
|H_{1}(\gamma^{*})\backslash O(\gamma)|+|H_2(\gamma^{*})| = \sum_{p \in V(\gamma^{*})}k(p).
\end{equation}
For $p \in V(\gamma^{*})$, let $k_1(p)$ be the number of vertical steps after $p$ in $\gamma^{*}$ intersecting with $l(p)$ and let $k_2(p)$ be the number of horizontal steps before $p$ in $\gamma^{*}$ intersecting with $l(p)$. 
Then $k(p) = k_1(p)+k_2(p)$ and 
\begin{align*}
\sum_{p \in V(\gamma^{*})}k_1(p) &= |H_{1}(\gamma^{*})\backslash O(\gamma)|\\
\sum_{p \in V(\gamma^{*})}k_2(p) &= |H_{2}(\gamma^{*})|.
\end{align*}
The reason why the first equation holds is that each of the inserted horizontal steps is just after an outer vertex, and the reason why the second equation holds is that each vertical step in $\gamma^{*}$ is just before an outer vertex.
Therefore, we have (\ref{ess}). \qed
\end{mainprf}

\end{document}